\documentclass{ws-ijmpe}
\usepackage[super,compress]{cite}
\usepackage{multirow}
\usepackage{float}
\usepackage[colorlinks=true,citecolor=blue,urlcolor=blue,linkcolor=red]{hyperref}

\begin{document}

\catchline{}{}{}{}{}

\title{Constraining Neutron Star Properties and Dark Matter Admixture with the NITR-I Equation of State: Insights from Observations and Universal Relations}
\author{
Pinku Routaray$^{1}$,
H.C. Das$^{2}$,
Jeet Amrit Pattnaik$^{3}$,
Bharat Kumar $^{1}$}

\address{$^{1}$ Department of Physics \& Astronomy, National Institute of Technology, Rourkela 769008, India}
\address{$^{2}$ INFN Sezione di Catania, Dipartimento di Fisica, Via S. Sofia 64, 95123 Catania, Italy}
\address{$^{3}$ Department of Physics, Siksha $'O'$ Anusandhan, Deemed to be University, Bhubaneswar -751030, India}

\maketitle

\begin{history}
\end{history}

\begin{abstract}
 
A recent observational study has set a constraint on the maximum mass of neutron stars (NSs), specifically focusing on PSR J0952-0607 and the compact star remnant HESS J1731-347, particularly within the low-mass regime. In our recent study \cite{pinku_jcap_2023}, we developed an energy density functional named NITR, which successfully produced the mass limit of the aforementioned pulsar but did not fully meet other observational constraints, such as those from NICER+XMM and GW170817. In this study, we introduce a new EDF named ``NITR-I," which not only reproduces the mass limit of PSR J0952-0607 but also aligns its canonical radius with NICER+XMM data, and its canonical dimensionless tidal deformability is consistent with the GW170817 event, thereby demonstrating the robustness of our model. The low-mass constraint associated with HESS J1731-347 suggests various possible compositions for the NS. The NITR-I model alone does not satisfy the HESS J1731-347 constraint; thus, we explore the possibility of incorporating dark matter (DM) inside the NS to meet this constraint. This approach proves successful when a specific value of Fermi momentum is considered. We also examine the impact of DM with varying Fermi momentum on different NS properties, such as tidal deformability and non-radial $f$-mode oscillation, using various relativistic mean-field (RMF) models. For the NITR-I EOS, the $f$-mode frequency is about 2.15 kHz at 1.4 $M_\odot$ when $k_f^{\rm DM} = 0$ GeV, and it slightly increases to around 2.32 kHz with $k_f^{\rm DM} = 0.03$ GeV. This increase in frequency due to DM suggests a possible reduction in tidal deformability, indicating that neutron stars with higher DM content are less susceptible to deformation by tidal forces which could be detectable in gravitational wave signals from neutron star mergers. Finally, we explore various universal relations (URs) for DM-admixed NSs, such as the relation between compactness and tidal deformability, the $f$-mode frequency and tidal deformability, and estimate the canonical values corresponding to both compactness and $f$-mode frequency using the GW170817 data.

\end{abstract}


\ccode{PACS numbers: 21.65.+f;  26.60.+c; 95.35.+d; 04.30.-w}

\section{Introduction}
\label{intro}
In recent decades, significant progress has been made in understanding the properties of compact stars, including their mass, radius, and tidal deformability \cite{Abbott_2017, Abbott_2018, Miller_2019, Miller_2021}. This advancement has both enriched the exploration of highly dense neutron stars (NS) and constrained various models. One pivotal breakthrough came from the observation of gravitational waves (GW) during the merger of two NSs in the GW170817 event \cite{Abbott_2017, Abbott_2018}. This observation provided valuable insights into the masses of binary systems and their tidal deformabilities. Another binary coalescence event, GW190814 \cite{RAbbott_2020}, involving a black hole and a compact object with a mass range of $2.50-2.67 \ M_\odot$, has sparked a compelling debate about the nature of the secondary compact object—whether it is the lightest black hole or a massive NS.

In addition to GW observations, concurrent measurements of the mass and radius of PSR J0030+0451 have been made through pulse profile modeling of X-ray emissions from hot spots on the surface of the isolated NS observed by the Neutron Star Interior Composition Explorer (NICER) telescope \cite{Riley-nicer_2019, Miller-nicer_2019}. The discovery of the Galactic NS PSR J0952-0607, named the ``Black widow pulsar," was unveiled in the Milky Way disc. It has been distinguished as the fastest and heaviest of its kind, with a mass of $M=2.35\pm 0.17 \ M_\odot$ \cite{Romani_2022}. Recently, a central compact object within the supernova remnant HESS J1731-347 has been investigated through the analysis of its X-ray spectrum, revealing it to be a low-mass object \cite{HESS_2022}. Remarkably, its mass is determined to be $M=0.77^{+0.20}_{-0.17} \ M_\odot$, with a radius of $10.4^{+0.86}_{-0.78}$ km, which opens a debate about whether this object is either the lightest NS or a strange star with a more exotic equation of state (EOS).

Several studies have already hypothesized about the nature of the HESS object. For example, Horvath {\it et al.} \cite{Horvath-hess_2023} employed the quark model and reached the conclusion that the object might be a light strange star. The notion of a strange nature for HESS J1731-347 was also put forth in Refs. \cite{clemente_2023hess, hcdas_2023hess, rather_2023quark}. On the other hand, some works proposed that it might be the lightest NS. For instance, Kubis {\it et al.} incorporated the $\omega-\delta$ cross-coupling and established a new relativistic mean-field (RMF) model by varying the slope parameter, which satisfied the HESS constraints \cite{kurbis_2019_prcmodel,kubis_2023hess}. Additionally, Huang {\it et al.} \cite{huang_2023hess} explored the possibility of the lightest NS by considering the RMF model with tensor coupling. A comprehensive study by Sagun {\it et al.} \cite{sagun_2023hess} investigated the nature of the central compact object (CCO) in HESS as baryonic, strange, or even admixed with dark matter (DM), particularly using a two-fluid approach. Hence, there is still a window to explore the nature of the HESS object.

In our recent study \cite{pinku_jcap_2023}, we developed an energy density functional (EDF) named ``NITR" within the formalism of RMF approximation. The NITR model was able to produce a maximum mass of $2.355 \ M_\odot$ in the light of PSR J0952-0607 \cite{Romani_2022}. However, it could not fully satisfy other observational constraints, such as the canonical radius measured by NICER+XMM data \cite{Miller_2021} and the canonical tidal deformability obtained from GW170817 data \cite{Abbott_2018}. Additionally, it was unable to reproduce the mass bounds suggested by HESS J1731-347. To address these limitations, we have devised a new parameter set, called NITR-I, for modeling an EDF based on the RMF approximation. The NITR-I model successfully produces the maximum mass consistent with PSR J0952-0607, meets the canonical radius constraint by NICER+XMM data, and aligns with the tidal deformability value from the GW170817 event, demonstrating the robustness of this model. However, the current model also fails to meet the low mass requirement predicted by HESS J1731-347 \cite{HESS_2022}. Therefore, to examine the low mass nature of the HESS object, we explore a potential approach by incorporating dark matter (DM) within the neutron star to meet this constraint. This approach proves successful when a specific value of Fermi momentum is considered, as discussed in detail in Sec.\ref{sec:RD}.

Compact stars, such as white dwarfs, NS, and strange stars, could efficiently capture DM particles due to their large baryonic density and very high gravitational potential \cite{McCullough_WD-DM_2010, Goldman_1989, Sandin_2009, Kouvaris_2008, Kouvaris_2011, Grigorious_2017, Lopes-harish_2023}. There are different ways of capturing DM inside the star, which may occur throughout their lifetime. The captured DM particles can modify the properties of the compact star depending on their interactions with nucleons and the mass of DM. Hence, here we examine the properties of DM admixed NS (DMANS) for the new NITR-I model along with other well-known RMF models. The existence of DM affects various NS properties like mass and radius; therefore, we should also look into its effects on the tidal properties and oscillation of the star. The oscillation is an essential tool that can be utilized to investigate the internal structure of the star, and it can be either radial \cite{Kokkotas_2001,souhardya_2023, pinku_prd_2023, Ishfaq_rad-delta_2023} or non-radial \cite{Bikram_2021-fmode,athul_2022, sailesh_2024, Probit_2024}, still, we have difficulties in detecting the NS oscillation modes. Therefore, studying the universal relation (UR) is a unique way to estimate those quantities which are not directly observed \cite{Kent_yagi_2013, Kent_yagi_2015, Gupta_2017, Yeung_2021, Harish__I_LOVE_C_2022, sailesh_2024}. Hence, in this study we investigate various URs for DMANS such as $C-\Lambda$ and also examine the $C-f$, $f-\Lambda$ URs within the relativistic Cowling framework \cite{sotani_2011-cowling, Flores_2014-cowling}. Also with the help of GW170817 tidal deformability data, we try to estimate the canonical compactness and the $f$-mode frequency of the star for both with and without DM.
\section{Formalism}
\label{sec:form}

\subsection{Nucleon EOS within RMF formalism}
\label{NITR-I model}

Here, we present the formalism that was followed to develop our new model, NITR-I. To achieve this, we consider the following Lagrangian density \cite{FURNSTAHL_1996, Frun_1997, singh_2014, Kumar_2017, Kumar_2018, Parmar_2022_1},

\begin{align}
{\cal L}_{\rm nucl.} & =   \sum_{\alpha=p,n} \bar\psi_{\alpha}
\Bigg\{\gamma_{\mu}\bigg(i\partial^{\mu}-g_{\omega}\omega^{\mu}-\frac{1}{2}g_{\rho}\vec{\tau}_{\alpha}\!\cdot\!\vec{\rho}^{\,\mu}\bigg) -\bigg(M_N - g_{\sigma}\sigma\bigg)\Bigg\} \psi_{\alpha} +\frac{1}{2}\partial^{\mu}\sigma\,\partial_{\mu}\sigma
\nonumber \\
&
-\frac{1}{2}m_{\sigma}^{2}\sigma^2+\frac{\zeta_0}{4!}g_\omega^2(\omega^{\mu}\omega_{\mu})^2-\frac{\kappa_3}{3!}\frac{g_{\sigma}m_{\sigma}^2\sigma^3}{M_N}-\frac{\kappa_4}{4!}\frac{g_{\sigma}^2m_{\sigma}^2\sigma^4}{M_N^2}+\frac{1}{2}m_{\omega}^{2}\omega^{\mu}\omega_{\mu}
\nonumber \\
&
-\frac{1}{4}W^{\mu\nu}W_{\mu\nu}
+\frac{1}{2}m_{\rho}^{2}\bigg(\vec\rho^{\mu}\!\cdot\!\vec\rho_{\mu}\bigg)  -\frac{1}{4}\vec R^{\mu\nu}\!\cdot\!\vec R_{\mu\nu} - \Lambda_{\omega}g_{\omega}^2g_{\rho}^2\big(\omega^{\mu}\omega_{\mu}\big)\big(\vec\rho^{\,\mu}\!\cdot\!\vec\rho_{\mu}\big)\, .
\end{align}

The presented Lagrangian density includes the nucleon fields for both protons ($p$) and neutrons ($n$), denoted as $\psi_\alpha$ ($\alpha = p,n$). Additionally, it incorporates the exchange meson fields $\sigma$, $\omega$, and $\rho$, along with their corresponding masses ($m_\sigma$, $m_\omega$, and $m_\rho$) and coupling constants ($g_\sigma$, $g_\omega$, and $g_\rho$). The terms $k_3$ and $k_4$ represent the self-interaction of the $\sigma$ meson, while $\zeta_0$ represents the self-interaction of the $\omega$ meson. The parameter $\Lambda_\omega$ defines the cross-coupling between the $\omega$ and $\rho$ mesons. Furthermore, various tensor fields, such as $W^{\mu\nu}$ and $\vec{R}^{\mu\nu}$, are also included.

The above Lagrangian density can be solved in the mean-field approach using both beta equilibrium and charge neutrality conditions, and one can calculate the energy density (${\cal{E}}_{\rm NS}$) and pressure ($P_{\rm NS}$) for the NS \cite{Kumar_2018}.
\begin{align}
{\cal{E}}_{\rm NS}&=\frac{2}{(2\pi)^{3}}\int d^{3}k E_{i}^\ast (k)+\rho  W+
\frac{ m_{s}^2\Phi^{2}}{g_{s}^2}\Bigg(\frac{1}{2}+\frac{\kappa_{3}}{3!}
\frac{\Phi }{M_N}+ \frac{\kappa_4}{4!}\frac{\Phi^2}{M_N^2}\Bigg) -\frac{1}{2}m_{\omega}^2\frac{W^{2}}{g_{\omega}^2}
\nonumber\\
&
-\frac{1}{4!}\frac{\zeta_{0}W^{4}} {g_{\omega}^2}+\frac{1}{2}\rho_{3} R -\frac{1}{2}\frac{m_{\rho}^2}{g_{\rho}^2}R^{2} -\Lambda_{\omega}  (R^{2}\times W^{2}) + \mathcal{E}_l \ \,
\label{enrns}
\end{align}
\begin{align}
P_{\rm NS} & = \frac{2}{3 (2\pi)^{3}}\int d^{3}k \frac{k^2}{E_{i}^\ast (k)}- \frac{m_{s}^2\Phi^{2}}{g_{s}^2}\Bigg(\frac{1}{2}+\frac{\kappa_{3}}{3!}\frac{\Phi }{M_N}+ \frac{\kappa_4}{4!}\frac{\Phi^2}{M_N^2}  \Bigg) +\frac{1}{2}m_{\omega}^2\frac{W^{2}}{g_{\omega}^2} 
\nonumber\\ 
& 
+\frac{1}{4!}\frac{\zeta_{0}W^{4}}{g_{\omega}^2} +\frac{1}{2} \frac{m_{\rho}^2}{g_{\rho}^2}R^{2} +\Lambda_{\omega} (R^{2}\times W^{2}) + P_l \ \,
\label{presns}
\end{align}
where $\Phi$, $W$ and $R$ are the fields associated with $\sigma$, $\omega$ and $\rho$ mesons respectively. $M_N$ is the mass of the nucleon. The ${\cal E}_l$ and  $P_l$ are the energy density and pressure of the leptons ($e^-, \mu^-$). 

The NITR-I model employed in this study was successfully calibrated using the simulation annealing method, allowing for precise determination of the various parameters of the RMF model, as mentioned above. The detailed parameterization procedure has been extensively described in Refs. \cite{BKAgrawal_2005, BKAgrawal_2006, Kumar_2017, Kumar_2018}. To achieve this calibration, coupling constants and nuclear matter properties were carefully fitted by employing experimental data on binding energy and charge radii for $^{16}$O, $^{40}$Ca, $^{48}$Ca, $^{68}$Ni, $^{90}$Zr, $^{100}$Sn, $^{132}$Sn, and $^{208}$Pb nuclei. The resulting values of these calibrated parameters are summarized in Table \ref{tab:parameter}. This rigorous calibration process ensures that the NITR-I model accurately represents the properties of neutron stars, providing a reliable and consistent framework for this investigation.
\begin{table*}
\centering
\caption{Different parameters of the ``NITR-I" model and its nuclear matter properties are given. Mass of the $\sigma$, $\omega$, and $\rho$ mesons are in MeV, and the nucleon mass is taken as 939 MeV.}
\renewcommand{\tabcolsep}{0.01cm}
\renewcommand{\arraystretch}{1.0}
\scalebox{0.8}{
\begin{tabular}{ccccccccccc}
\hline \hline
Model & $m_\sigma$ &  $m_\omega$   & $m_\rho$   & $g_\sigma$   & $g_\omega$   &$g_\rho$   & $\kappa_3$   &$\kappa_4$  & $\zeta_0$   & $\Lambda_{\omega}$\\
\hline
 NITR-I & 470.00 & 782.5 & 763.0  & 8.729  & 11.172  & 9.461  & 2.729  & -10.207  & 0.159  & 0.029 \\ [0.5cm]
\hline
$\rho_{\rm sat.}$ (fm$^{-3}$) &  & $\mathcal{E}_{\rm sat.}$ (MeV) &  &$K_{\rm sat.}$ (MeV)&  &$J_{\rm sat.}$ (MeV) & &$L_{\rm sat.}$ (MeV) & & \\ \hline
0.151 &  &-16.337 &  &199.018 &  &30.937 & &61.826 & &  \\
\hline \hline
\end{tabular}}
\label{tab:parameter}
\end{table*}

The structure of the NS is dispersed from the crust to the core. In essence, compared to the crust, the core is distributed throughout the densest region of the NS and contributes the majority of the NS's mass. Since the crust acts as an intermediary between the core and the surface, all information originating from the core gets impacted. To provide consistency and a more precise assessment of the characteristics of the NSs, we need an EOS that is considered from core to crust region regarded as ``Unified EOS''. Therefore, the crust EOS for the NITR-I model is obtained by following the methodology used in Ref. \cite{Parmar_2022}. By patching the same crust EOS with the core part, we obtain the unified EOS for the NITR-I model, which is used to explore the properties of the NS.
\begin{figure*}
   \centering
   \includegraphics[width = 0.9\textwidth]{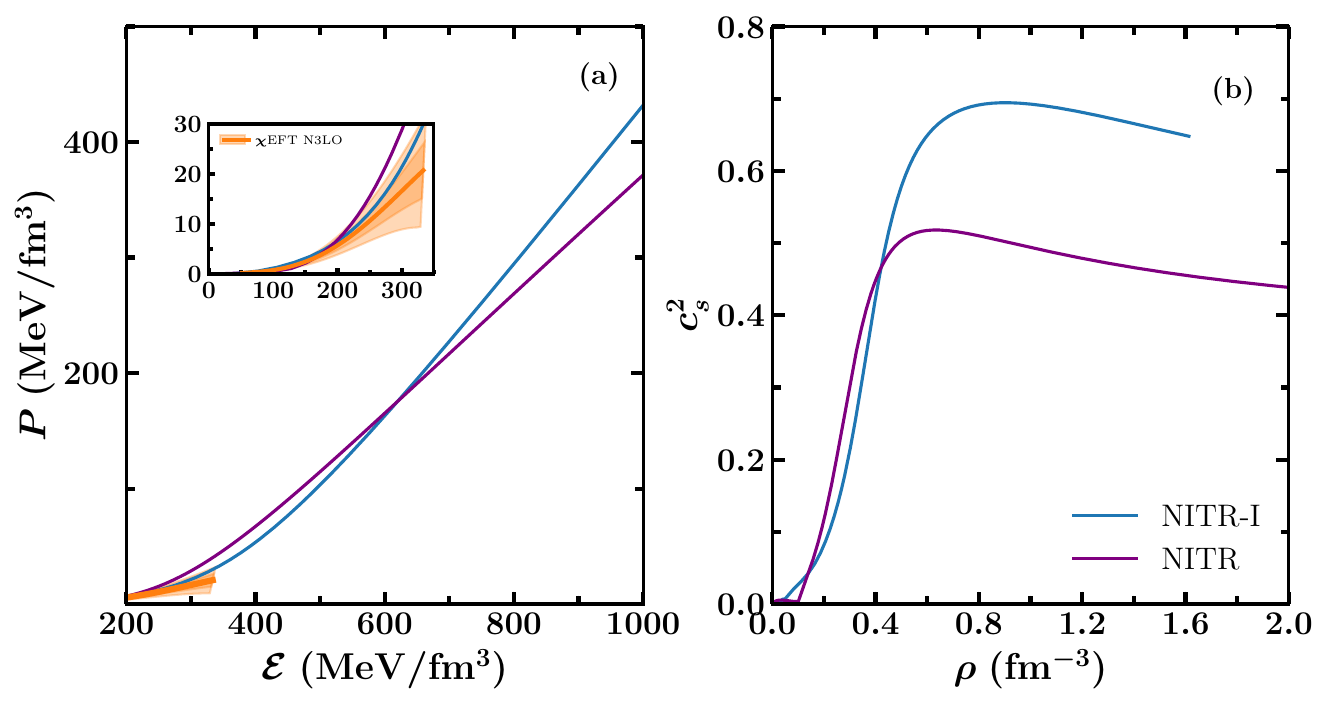}
    \caption{{\bf (a)} The EOS for the NITR-I model is displayed and compared alongside with NITR model. The low-density constraint by the chiral Effective Field Theory (EFT) is also shown \cite{Drischler_2021}. {\bf (b)} The sound speed ($C_s^2$) as a function of the baryon density.}
    \label{fig:NITR-I_EOS_Cs}
\end{figure*}

The unified EOSs for the NITR-I and NITR models are displayed in the left panel of Fig. \ref{fig:NITR-I_EOS_Cs}(a). The EOS for the NITR-I model demonstrates better agreement with the EFT constraints at low-density regions, as indicated by the overlapping regions with the orange EFT band. This consistency suggests that the NITR-I model is more robust and reliable in describing neutron-rich matter at lower densities compared to the NITR model. The difference between the models becomes more pronounced at higher energy densities, where the NITR-I model predicts a stiffer EOS, leading to higher pressures. The right panel, Fig. \ref{fig:NITR-I_EOS_Cs}(b), illustrates the behavior of the squared speed of sound $ C_s^2 $ as a function of baryon density $ \rho $. For both models, the speed of sound initially rises with increasing density but exhibits distinct characteristics at higher densities. The NITR-I model shows a peak in $ C_s^2 $ at a baryon density around $ \rho \approx 0.8 \, \text{fm}^{-3} $, after which it slightly decreases, while the NITR model demonstrates a more gradual rise followed by a flattening. Notably, both models satisfy the causality condition $ C_s^2 < 1 $ throughout the range, ensuring that the sound speed remains below the speed of light, thus validating the physical realism of these models.

\begin{figure*}
   \centering
    \includegraphics[width = 0.9\textwidth]{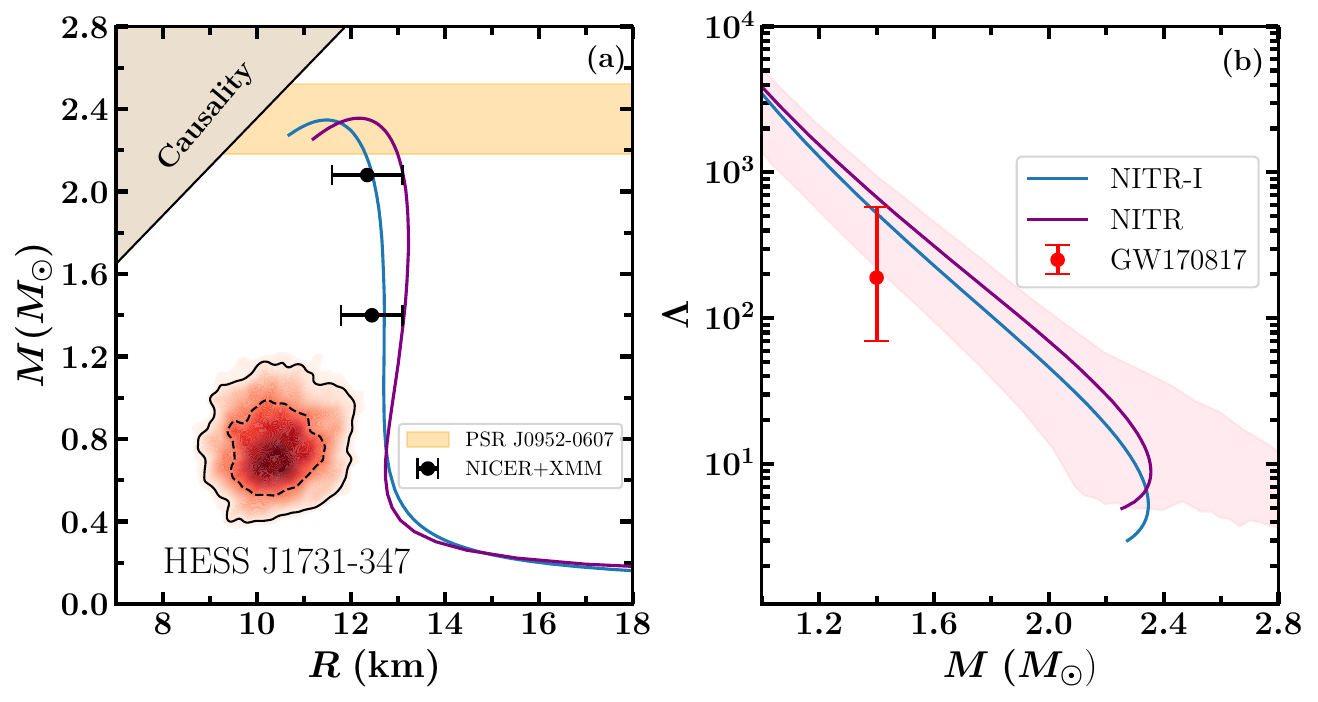}
    \caption{\textbf{(a)} The mass-radius profiles of the NS are shown forNITR-I and NITR EOSs with various observational constraints such as recently observed heaviest pulsar PSR J0952-0607 \cite{Romani_2022}, HESS J1731-347 \cite{HESS_2022} and NICER+XMM data \cite{Miller_2021}. {\bf (b)} The tidal deformability is plotted as a function of the mass, including with the constraint from GW170817 event \cite{Abbott_2018}.}
    \label{fig:NITR-I_mr_tidal}
\end{figure*}
In Fig. \ref{fig:NITR-I_mr_tidal}, the mass-radius and tidal deformability profiles of NSs for the NITR-I and NITR models are compared against various observational constraints. The left panel (Fig. \ref{fig:NITR-I_mr_tidal}(a)) shows that the NITR-I model successfully reproduces the maximum mass of $ M = 2.34 \, M_\odot $, which is in excellent agreement with the recently observed heaviest pulsar PSR J0952-0607 , thus meeting the high-mass constraint. Both models also match the observational constraints for a canonical neutron star (with mass $ M = 1.4 \, M_\odot $) and the $ 2.08 \, M_\odot $ star as measured by NICER+XMM \cite{Miller_2021}. However, the NITR-I model performs slightly better in terms of reproducing the radius of the $ 2.08 \, M_\odot $ star, albeit both models fail to meet the mass constraint associated with HESS J1731-347 \cite{HESS_2022}, which suggests a significantly lower mass. The right panel (Fig. \ref{fig:NITR-I_mr_tidal}(b)) presents the tidal deformability $ \Lambda $ as a function of mass $ M $, with the observational constraint from the GW170817 event overlaid for reference \cite{Abbott_2018}. The NITR-I model exhibits a tidal deformability consistent with the GW170817 constraints, thereby validating its applicability within the context of binary neutron star mergers. In contrast, the NITR model produces values of $ \Lambda $ that fall outside the observed range, indicating less agreement with this crucial observational constraint. While both EOSs provide acceptable descriptions of some NS properties, they do not simultaneously satisfy all observational constraints, particularly the low-mass constraint from HESS J1731-347. This limitation suggests that relying solely on nucleonic degrees of freedom may be insufficient for a complete description of neutron star properties, particularly within the NITR-I model. As summarized in Table \ref{tab:NITR-I_results}, key NS properties such as maximum mass, corresponding radius, and tidal deformability are presented for both the NITR-I and NITR models. To address discrepancies, the inclusion of DM within the neutron star is proposed as a plausible solution. The detailed exploration of this scenario is presented in the subsequent sections.

\begin{table}
\centering
\caption{The maximum mass ($M_{\rm max}$) and its corresponding radius ($R_{\rm max}$) for NITR-I model and also the canonical radius ($R_{1.4}$) and tidal deformability ($\Lambda_{1.4}$) are given. The results for the NITR model are given for comparison.}
\renewcommand{\tabcolsep}{0.2cm}
\renewcommand{\arraystretch}{1.5}
\begin{tabular}{ccccc}
\hline \hline
Model & $M_{\rm max} \ (M_\odot)$ &  $R_{\rm max}$(km)   & $R_{1.4}$(km)   & $\Lambda_{1.4}$ \\
\hline
 NITR-I & 2.34 & 11.46  & 12.71  & 526.96 \\
 NITR & 2.35 & 12.19  & 13.13  & 682.84 \\
\hline \hline
\end{tabular}
\label{tab:NITR-I_results}
\end{table}
\subsection{Dark Matter EOS using single-fluid approach}
The exploration of the true nature of DM remains an ongoing challenge in the field of cosmology. Over the course of many years since the Big Bang, numerous DM candidates have been proposed, while weakly interacting massive particles (WIMPs) have gained significant popularity due to their thermal relic nature. Our study focuses on a specific class of WIMPs, known as non-annihilating WIMPs or neutralinos, which we consider a potential candidate for DM. These neutralinos are hypothesized to have already accreted within the NS due to their high baryon density and enormous gravitational potential \cite{Kouvaris_2011, Goldman_1989, arpan_2019}.

To elucidate the nature of DMANS, we construct a Lagrangian density that considers the coupling between baryons and DM through Higgs exchange. The formulation is motivated by previous works \cite{Grigorious_2017, arpan_2019, harishmnras_2020, pinku_prd_2023, pinku_jcap_2023, pinku_mnras_2023} and is as follows:
\begin{align}
{\cal{L}}_{\rm DM} & =  \bar \chi \left[ i \gamma^\mu \partial_\mu - M_\chi + y h \right] \chi +  \frac{1}{2}\partial_\mu h \partial^\mu h - \frac{1}{2} M_h^2 h^2 + \sum_{\alpha=p,n} f \frac{M_N}{v} \bar \psi_\alpha h \psi_\alpha \, .
\label{eq:LDM}
\end{align}
The symbols $\psi_\alpha$ and $\chi$ represent the nucleon and DM wave functions, respectively. The Higgs field is denoted by the symbol $h$. The masses $M_\chi$ and $M_h$ represent the neutralino and Higgs mass, respectively. The coupling constants between the DM and SM Higgs Bosons are denoted as $y$. Consequently, $f M_N/v$ denotes the nucleon-Higgs field effective Yukawa coupling with proton-Higgs form factor $f$ and vacuum expectation of Higgs $v$.The values for all parameters used in this study are given in Refs. \cite{harishmnras_2020, Das_fmode_2021}.

The above Lagrangian can be solved to obtain the EOS for the DM and is given by \cite{Grigorious_2017,harishmnras_2020,pinku_prd_2023},
\begin{eqnarray}
{\cal{E}}_{\rm DM} = \frac{2}{(2\pi)^{3}}\int_0^{k_f^{\rm DM}} d^{3}k \sqrt{k^2 + {M_\chi^\star}^2 } + \frac{1}{2}M_h^2 h_0^2 \, ,
\label{eq:edm}
\end{eqnarray}
\begin{eqnarray}
P_{\rm DM} = \frac{2}{3(2\pi)^{3}}\int_0^{k_f^{\rm DM}} \frac{d^{3}k \hspace{1mm}k^2} {\sqrt{k^2 + {M_\chi^\star}^2}} - \frac{1}{2}M_h^2 h_0^2 \, ,
\label{eq:pres}
\end{eqnarray} 
where $M_\chi^\star$ represents the effective mass of the DM, and $k_f^{\rm DM}$ is the DM Fermi momentum. 

Hence, in the context of DMANS, the expressions for the total energy density and pressure can be formulated as follows. 
\begin{eqnarray}
\mathcal{E}_{\rm DMANS}={\cal{E}}_{\rm NS}+ {\cal{E}}_{\rm DM} \, ,
\nonumber
\\
{\rm and}
\hspace{1cm}
P_{\rm DMANS}=P_{\rm NS} + P_{\rm DM} \, .
\label{eq:EOS_total}
\end{eqnarray}

Once the EOS is determined as an input, the Tolman–Oppenheimer–Volkoff (TOV) equation can be solved to compute the mass and radius of the DMANS \cite{Tolman_1939, Oppenheimer_1939, pinku_prd_2023}. The calculation for the dimensionless tidal deformability ($\Lambda$) \cite{Hinderer_2008, Kumartidal_2017} and non-radial $f$-mode oscillation \cite{athul_2022,Das_fmode_2021, Bikram_2021-fmode,Probit_2024} can be done by solving their respective differential equations along with the TOV equation as done in Ref. \cite{Das_fmode_2021}.
NS properties with admixed DM
\label{sec:RD}

In this section, we investigate various macroscopic properties of the NS and derive the effects of DM on it. Firstly, we focus on determining the nature of HESS J1731-347, then we investigate the tidal deformability and $f$-mode frequency. In order to explain the robustness of this study, we choose our new NITR-I model along with varieties of well-known RMF models, such as FSU2 \cite{Chen_2014}, FSUGarnet \cite{Chen_2015}, G1 \cite{Frun_1997}, IOPB-I \cite{Kumar_2018}, NITR \cite{pinku_jcap_2023}, and TM1 \cite{Sugahara_1994}. All the models are able to satisfy the $2 M_\odot$ constraint.

\subsection{HESS J1731-347 might be a DMANS}
As discussed in the earlier section on the NITR-I model, the newly developed NITR-I EOS satisfies most of the observational bounds but fails to achieve the mass limit provided by HESS J1731-347. To address this issue, we explore the impact of DM on the mass-radius relation by varying the DM content within the star. The inclusion of DM softens the EOS of the neutron star, which subsequently reduces its macroscopic properties such as mass and radius. 

\begin{figure*}
   \centering
   \includegraphics[width = 1.0\textwidth]{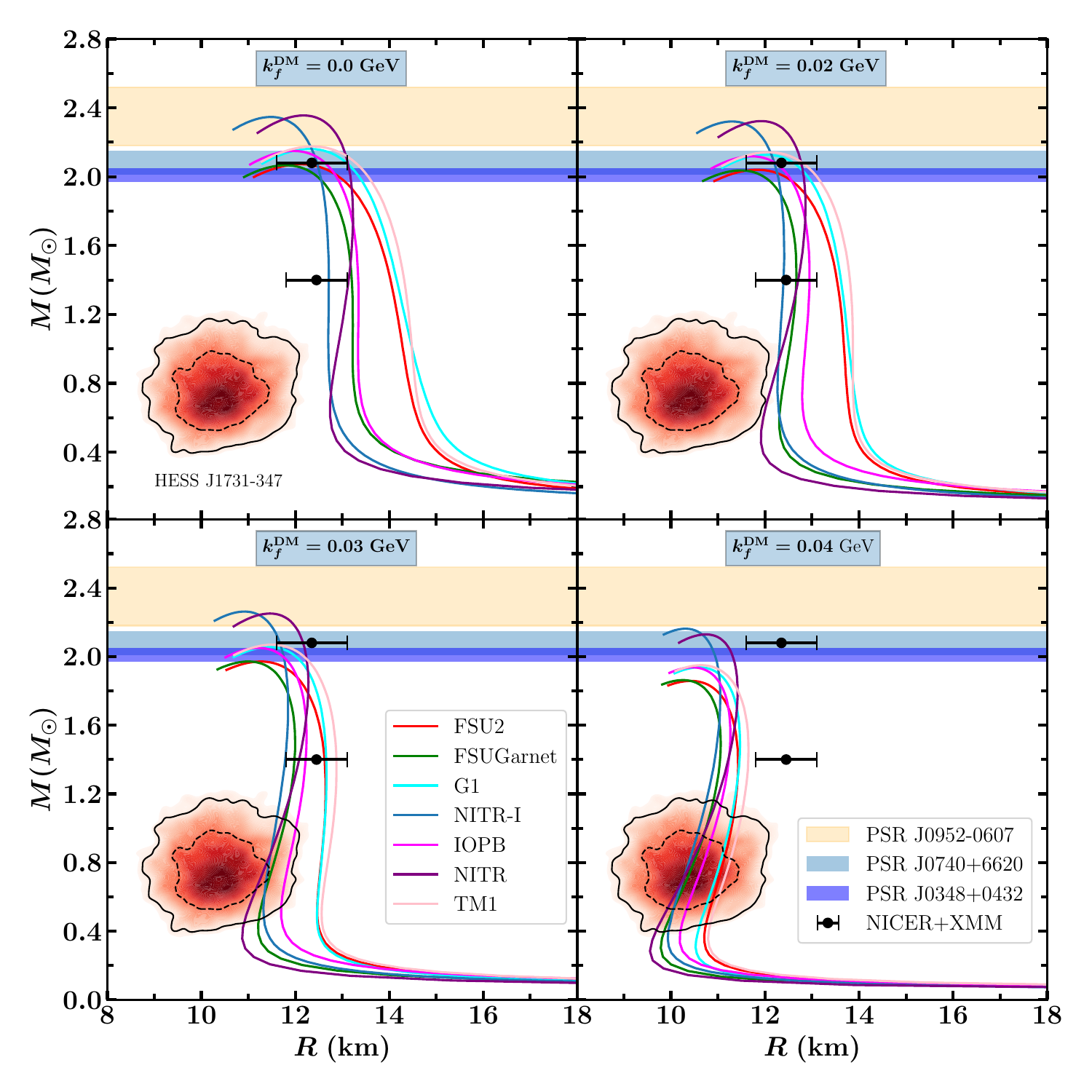}
    \caption{The mass-radius profiles are shown with and without DM by varying $k_f^{\rm DM}$. The color bands represent the mass range of different pulsars such as PSR J0952-0607 \cite{Romani_2022}, PSR J0740+6620 \cite{Cromartie_2020}. The black error bars represent the radii limits for the canonical star and $2.08 \ M_\odot$ star reported by NICER+XMM data \cite{Miller_2021}. The 1$\sigma$ and 2$\sigma$ contours of HESS J1731-347 \cite{HESS_2022} are represented by dashed lines and solid lines, respectively, in the contour.}
    \label{fig:mr_DM}
\end{figure*}

In Fig. \ref{fig:mr_DM}, the mass-radius profiles are depicted both with and without DM, alongside various observational constraints. Without DM ($ k_f^{\rm DM} = 0.0 \, \text{GeV} $), none of the RMF models can satisfy the HESS data. However, with a DM momentum of $ k_f^{\rm DM} = 0.03 \, \text{GeV} $, the NITR-I model, along with a few other RMF models, is able to satisfy all the observational constraints simultaneously, including the HESS data. At a higher DM momentum ($ k_f^{\rm DM} = 0.04 \, \text{GeV} $), all RMF models reproduce the HESS bounds but deviate from the NICER+XMM data. These findings suggest that DM plays a significant role in reconciling the observational constraints, allowing us to estimate the amount of DM present inside neutron stars based on these constraints. Additionally, this study indicates that with a certain amount of DM, the NICER+XMM and HESS data can be simultaneously satisfied. This also aligns with the findings of Sagun et al. \cite{sagun_2023hess} , where they suggest that the HESS J1731-347 object might be a DMANS based on a two-fluid approach .

\subsection{Tidal deformability and Non-radial Oscillation of DMANS}
\begin{figure*}
   \centering
   \includegraphics[width = 1.0\textwidth]{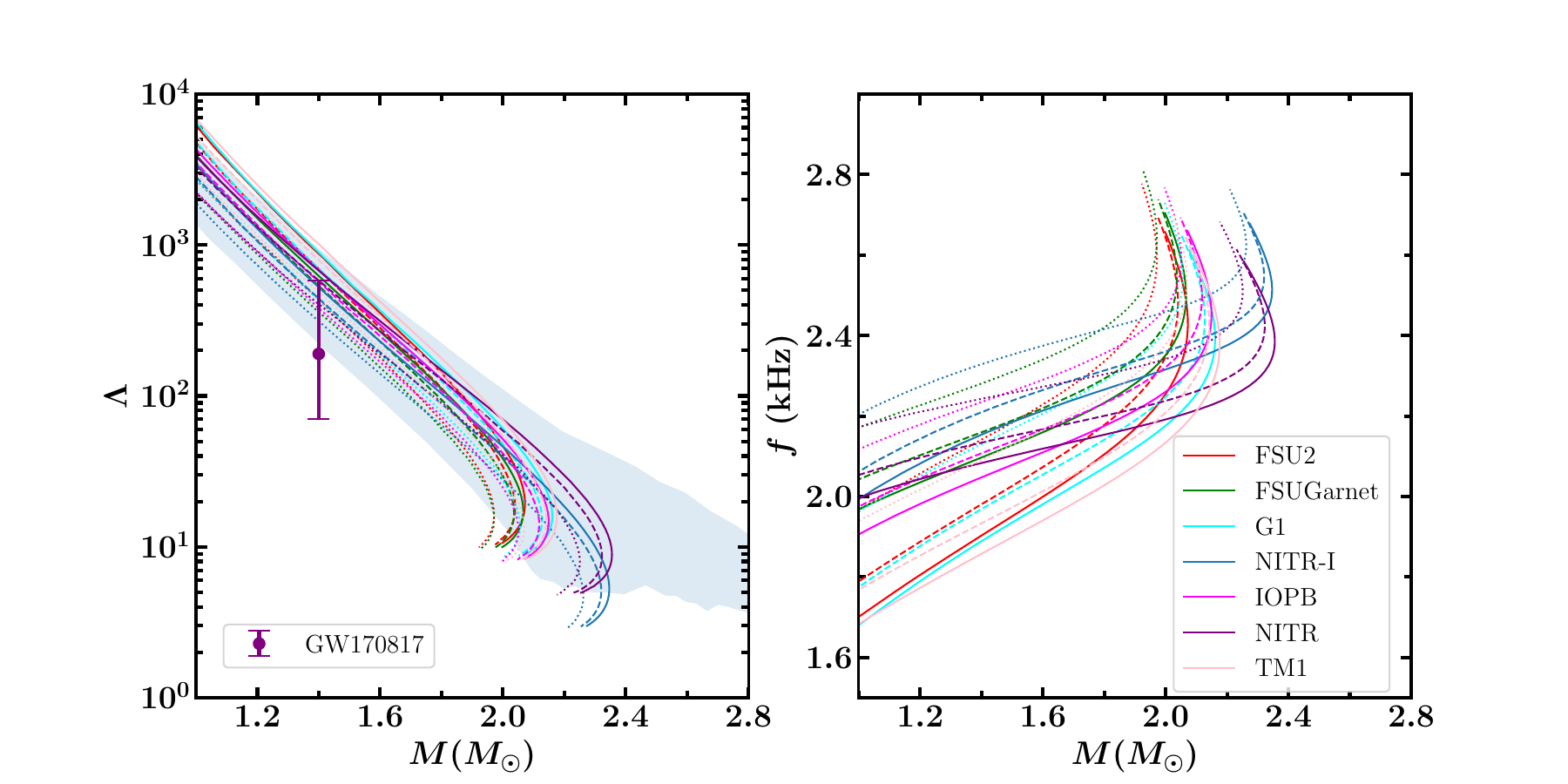}
    \caption{{\it Left:} Variation of the tidal deformability with mass for both with and without DM is shown, and the observational constraint for GW170817 \cite{Abbott_2018} is also imposed. The blue-shaded region corresponds to the constraints derived from the GW170817 measurements \cite{Abbott_2018}. {\it Right:} The non-radial $f$-mode frequency is shown as a function of the mass of the NS with different DM configurations. For both panels, the DM configurations are represented as $k_f^{\rm DM}=0.0 \ {\rm (solid \ line)}, \ 0.02 \ {\rm (dashed \ line) \ {\rm and} \ 0.03 \ {\rm GeV} \ {\rm (dotted \ line)}}$.}
    \label{fig:f-mode_tidal_DM}
\end{figure*}

As discussed earlier, the properties of NSs such as mass and radius are significantly impacted by the presence of dark matter. This influence extends to the tidal deformability and $f$-mode frequency of the star. In this subsection, we calculate the tidal deformability and $f$-mode frequency of DMANS for the considered EOSs, as shown in Fig. \ref{fig:f-mode_tidal_DM}.

The left panel of Fig. \ref{fig:f-mode_tidal_DM} presents the calculated tidal deformability values for both baryonic-only and DM-admixed neutron stars. The observational constraint from the GW170817 \cite{Abbott_2018} event is also shown, with the blue-shaded region corresponding to the constraints derived from this gravitational wave measurement . From the figure, it is evident that the tidal deformability decreases with the inclusion of DM. This decrease indicates that the presence of DM within the neutron star leads to a more compact structure, making the star less susceptible to deformation due to tidal forces. This finding aligns with earlier observations, reinforcing the notion that DM plays a critical role in altering the structural properties of neutron stars.

The right panel of Fig. \ref{fig:f-mode_tidal_DM} shows the variation of non-radial $f$-mode oscillation frequencies as a function of neutron star mass for various nuclear matter models, under different DM admixtures characterized by Fermi momentum $k_f^{\rm DM}$. The figure shows that the $f$-mode frequency generally increases with the mass of the neutron star across all equations of state. Moreover, for a fixed mass, the frequency rises with increasing $k_f^{\rm DM}$, indicating that a higher DM content results in denser configurations and consequently higher oscillation frequencies. Specifically, for the NITR-I EOS, the $f$-mode frequency is approximately 2.15 kHz at 1.4 $M_\odot$ with $k_f^{\rm DM} = 0$ GeV and increases slightly to around 2.32 kHz when $k_f^{\rm DM} = 0.03$ GeV. At the maximum mass, which is about 2.35 $M_\odot$ for the NITR-I EoS, the frequency is around 2.51 kHz for $k_f^{\rm DM} = 0$ GeV, and it rises to nearly 2.61 kHz for $k_f^{\rm DM} = 0.03$ GeV. Similar trends are observed in other EoSs, where the effect of DM becomes more pronounced at higher masses, resulting in a significant divergence of frequencies based on DM content. This increase in $f$-mode frequency with DM content suggests that neutron stars with more DM oscillate at higher frequencies, reflecting the stiffening of the star's core due to DM presence and leading to a corresponding increase in average density. Additionally, this rise in frequency with DM implies a potential decrease in tidal deformability, making neutron stars with higher DM content less deformable under tidal forces, which could have observable consequences in gravitational wave signals from neutron star mergers. The relationship between $f$-mode frequency and average density has been well established in previous studies \cite{Andersson-Kokkotas_GW_1996, Bikram_2021-fmode, Das_fmode_2021}. Moreover, universal relations linking tidal deformability and $f$-mode frequency across various EOSs have been identified  \cite{sotani-kumar-UR_2021,sailesh_2024}. In the following subsections, we will explore these relationships further and attempt to estimate the $f$-mode frequency using observational data.

\subsection{Relation between $f$-mode frequency and average density for DMANS}

In earlier works, the non-radial $f$-mode frequency has been derived as a function of the average density ($\bar{\rho}=\sqrt{\bar{M}/\bar{R}^3}$) of the star. In Ref. \cite{Andersson-Kokkotas_GW_1996}, Andersson and Kokkotas (AK) first presented an empirical relationship between these two variables using a limited number of polytropic EOSs. They later refined this relationship in Ref. \cite{Andersson-Kokkotas_asteroseismology_1998} by incorporating a broader range of realistic EOSs. In Ref. \cite{Das_fmode_2021}, Das (HCD) explored this relationship specifically for DM-admixed hyperonic stars, following Pradhan and Chatterjee's (PC) recent derivation for hyperonic stars in Ref. \cite{Bikram_2021-fmode}. Similarly, in our study, we have investigated the variation of the $f$-mode frequency with the average density for DMANS, extending the analysis to assess the influence of dark matter on this relationship.

\begin{figure*}
   \centering
   \includegraphics[width = 1.0\textwidth]{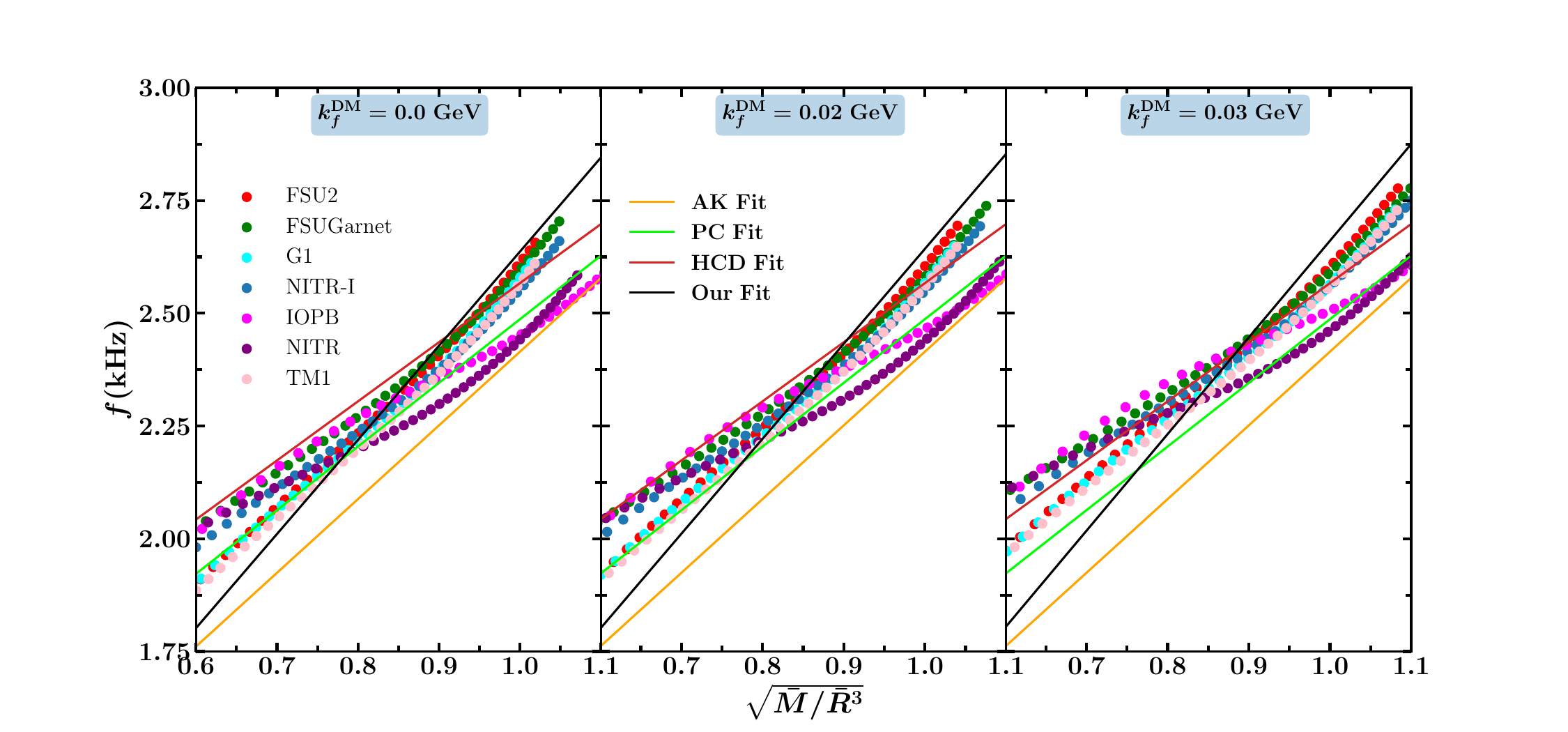}
    \caption{The $f$-mode frequency is fitted with the average density of the star using the Eq. \ref{eq:fit_f-n}. The best fit for our results with various DM configurations is represented by the black solid line in each panel, which we also compare with the previous works, including AK Fit \cite{Andersson-Kokkotas_GW_1996}, PC Fit \cite{Bikram_2021-fmode}, and HCD Fit \cite{Das_fmode_2021}.}
    \label{fig:avd_den-f_mode}
\end{figure*}

The empirical relation connecting $f$-mode frequency and the average density of the star can be written as follows,
\begin{eqnarray}
    f = a + b \sqrt{\frac{\bar{M}}{\bar{R}^3}}.
    \label{eq:fit_f-n}
\end{eqnarray}
Where $\bar{M}=\frac{M}{1.4 \ M_\odot}$ and $\bar{R}=\frac{R}{1.4 \ {\rm km}}$. $a$ and $b$ are fitting coefficients.

In Fig. \ref{fig:avd_den-f_mode}, we present the fitting of the $f$-mode frequency with the average density of neutron stars under various DM configurations, using the fitting equation provided in the text. The left panel of Fig. \ref{fig:avd_den-f_mode} displays the fitting relation for neutron stars without DM, based on several RMF EOSs. The calculated best fit is represented by the black solid line, and our results are compared with earlier fits as summarized in Table \ref{tab:avg_den-f_mode}. Notably, our model predicts a higher amplitude for the $f$-mode frequency compared to these earlier fits, indicating a more significant $f$-mode response for the EOSs considered in this study. In the middle panel, we incorporate DM with $k_f^{\rm DM}=0.02$ GeV and calculate the $f$-mode frequency for DM-admixed neutron stars. The best fit obtained using the empirical formula is again depicted by the black solid line. In comparison to the scenario without DM, we observe an increase in the amplitude of the $f$-mode frequency with rising average density. This suggests that the presence of DM increases the stiffness of the neutron star, leading to a higher oscillation frequency. The right panel explores the case where $k_f^{\rm DM} = 0.03$ GeV. Similar to the middle panel, the $f$-mode frequency continues to rise with increasing average density. This consistent behavior across panels underscores the trend that higher DM content results in a more pronounced increase in the $f$-mode frequency.

\begin{table}
\centering
\caption{Utilizing the empirical relation $f = a + b \sqrt{\frac{\bar{M}}{\bar{R}^3}}$, the corresponding obtained coefficients a and b are depicted with and without DM. Additionally, some earlier studies such as AK \cite{Andersson-Kokkotas_asteroseismology_1998}, PC \cite{Bikram_2021-fmode}, HCD \cite{Das_fmode_2021} fits are shown to compare with our results.}
\renewcommand{\tabcolsep}{0.05cm}
\renewcommand{\arraystretch}{1.2}
\begin{tabular}{ccc}
\hline \hline
Fittings & $a({\rm kHz})$ &  $b({\rm kHz})$\\
\hline
 AK Fit     &   0.78    &   1.635 \\
 PC Fit     &   1.075   &   1.412 \\
 HCD Fit    &   1.256   &   1.311 \\
 \hline
\hspace{5cm} This study &  &\\ \hline 
 $k_f^{\rm DM}=0.0 \ {\rm GeV}$ & 0.549 & 2.088 \\
 $k_f^{\rm DM}=0.02 \ {\rm GeV}$ & 0.537 & 2.106 \\
 $k_f^{\rm DM}=0.03 \ {\rm GeV}$ & 0.517 & 2.144 \\
\hline \hline
\end{tabular}
\label{tab:avg_den-f_mode}
\end{table}
\subsection{Universal Relations}
In astrophysics, universal relations (URs) are regarded as a crucial approach for theoretically determining the characteristics of stars. Since not all properties of a star can be determined from a single observation, URs allow us to derive unknown quantities as functions of known ones, as they are independent of the equation of state. In the literature, several URs have been developed to investigate the characteristics of neutron stars \cite{Kent_yagi_2013,Kent_yagi_2015,Gupta_2017,Yeung_2021,Chakrabarti_2014,Haskell_2013,Bandyopadhyay2018,Harish__I_LOVE_C_2022}.  In Refs. \cite{Prashant_2024}. In recent work \cite{Prashant_2024} , authours have explored the $C-\Lambda$ UR for DM-admixed NSs using a two-fluid approach, concluding that the presence of DM does not break the UR. However, there remains an open opportunity to analyze other URs for DM-admixed NSs, as DM can influence various macroscopic properties of these stars. Therefore, in this study, we present URs for $C-\Lambda$, $\Bar{\omega}-\Lambda$, and $C-\Bar{\omega}$ for both cases with and without DM. Here, $\Bar{\omega} = \omega M$ represents the normalized $f$-mode frequency, where $M$ is the mass and $\omega$ is the angular frequency.

\subsubsection{$C-\Lambda$ relation}

While analysing a gravitational wave signal generated during a binary coalescence, the $C-\Lambda$ relation may being a very helpful tool for obtaining insights on the NS EOS and this relation first identified by Maselli {\it et al.} \cite{Maselli_C-L_2013}. Here, we calculate the compactness and dimensionless tidal deformability for both with and without DM and the least-squares fitting is carried out utilising the approximate formula that follows, 

\begin{equation}\label{eq:C-L_fit}
C = \sum_{n=0}^{4} a_n \left( \log_{10}{\Lambda} \right)^n,
\end{equation}
where $a_n$ is the fitting coefficient that adjusts to the fourth order, and the corresponding coefficients are provided in Table \ref{tab:C-L-f_coef}.
 
\begin{figure*}
   \centering
   \includegraphics[width = 1.0\textwidth]{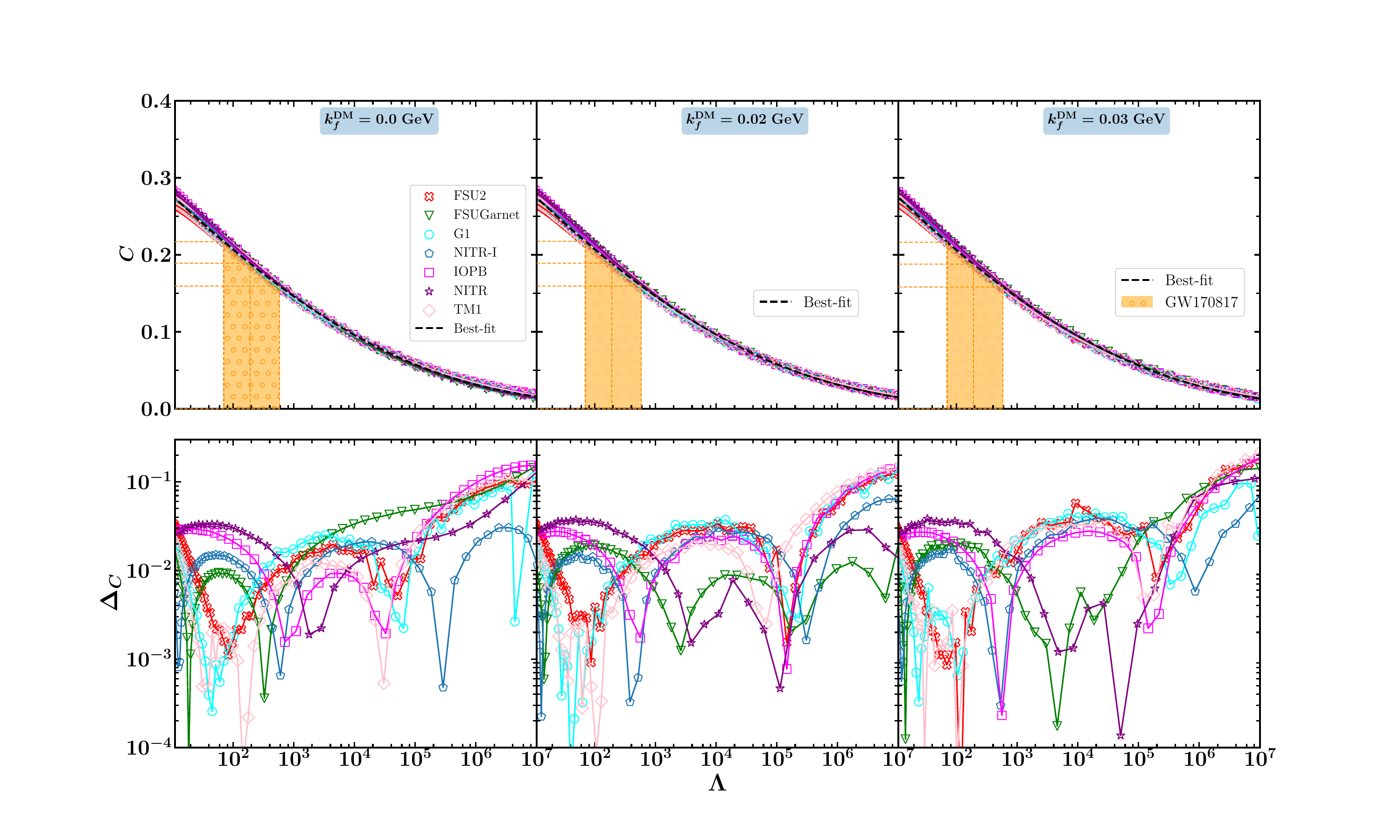}
    \caption{$\textbf{Top:}$ The UR relation between compactness and dimensionless tidal deformability is presented for both with and without DM. The best fit for all cases is depicted with a black dashed line. The canonical constraint from GW170817 \cite{Abbott_2018} is used to determine canonical compactness. $\textbf{Bottom:}$ The relative deviation for each case is shown in the lower panel.}
    \label{fig:C-L_UR}
\end{figure*}
In Fig. \ref{fig:C-L_UR} , we examine the UR between compactness ($C$) and dimensionless tidal deformability ($\Lambda$) across different DM configurations, represented by varying values of $k_f^{\rm DM}$. The upper panels display the $C-\Lambda$ relation, with the best fits depicted as black dashed lines. These fits correspond to cases with $k_f^{\rm DM} = 0.0$, $0.02$, and $0.03$ GeV, respectively. The figure clearly illustrates that the $C-\Lambda$ UR remains intact even with the inclusion of DM, indicating the robustness of this relationship. This is significant as it suggests that DM does not disrupt the underlying physics governing the compactness and tidal deformability of neutron stars. Furthermore, in the lower panels, the relative deviation in compactness ($\Delta C = \frac{|C - C_{\text{fit}}|}{C_{\text{fit}}}$) is plotted against $\Lambda$ for each case. Notably, the relative deviation remains below 1\% for most EOSs across all values of $\Lambda$, with a slight increase at high and low $\Lambda$ values, which reflects a strong agreement between our results and the best-fit line. Specifically, for the canonical NS model with $1.4 M_\odot$, the deviation in compactness is maintained under 3\% across all DM scenarios, both with and without DM. This consistency reinforces the reliability of the $C-\Lambda$ UR even in the presence of DM, allowing us to impose the canonical tidal constraint from GW170817 to estimate the canonical compactness ($C_{1.4}$) of the star. The influence of DM is further evidenced by the slight variations in $C_{1.4}$ values, as indicated in Table \ref{tab:C-L-f_coef}, highlighting the subtle yet measurable impact of DM on NS structure.

\subsubsection{$\Bar{\omega}-\Lambda$ relation}

\begin{figure*}
   \centering
   \includegraphics[width = 1.0\textwidth]{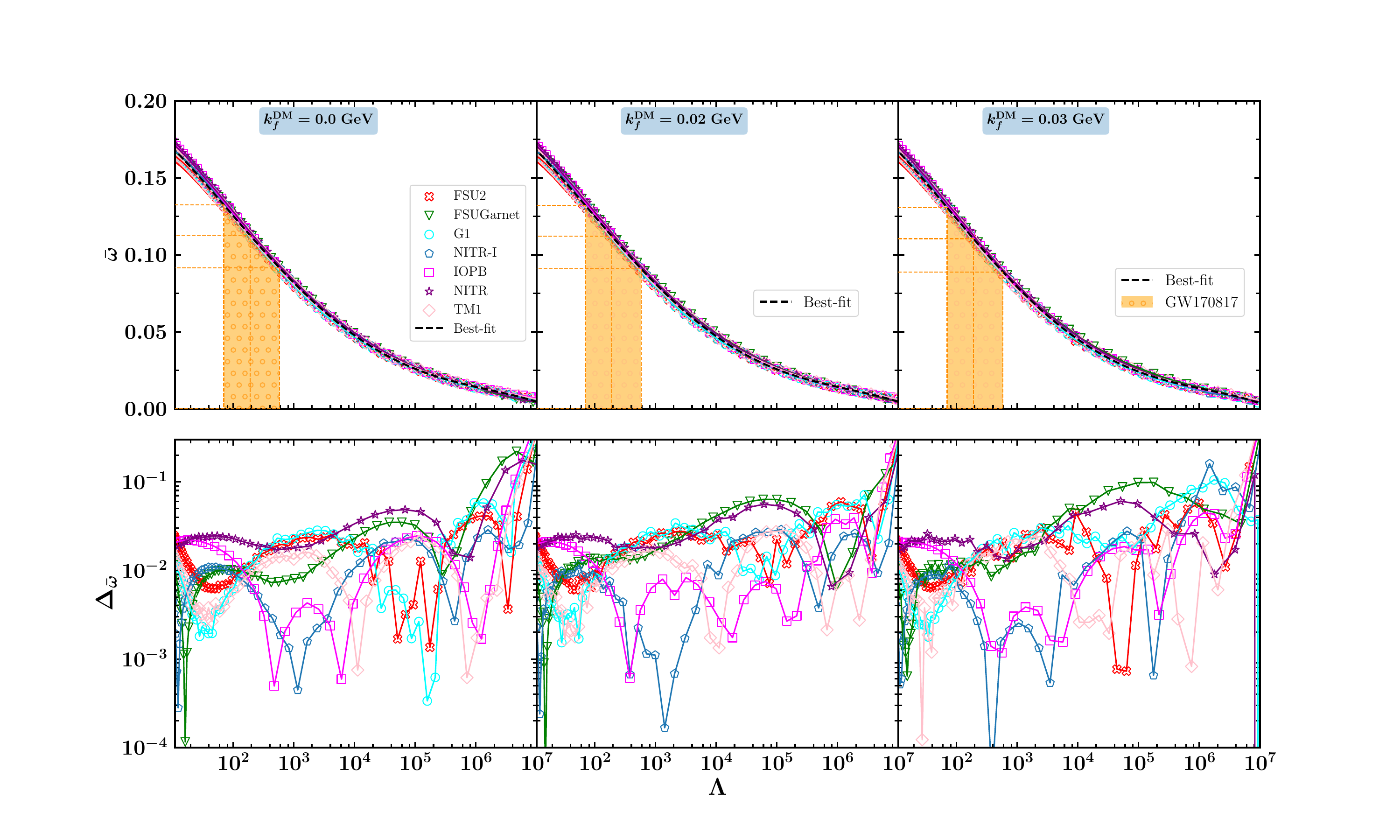}
    \caption{As similar to the Fig. \ref{fig:C-L_UR}, the UR is shown for $\Bar{\omega}-\Lambda$ relation.}
    \label{fig:F-L_UR}
\end{figure*}

The investigation of non-radial oscillation modes, such as the $f$-mode frequency in NSs, provides essential insights into their internal composition and the EOS governing their structure. Despite the challenges in detecting and analyzing these modes, URs serve as a valuable theoretical framework to explore these frequencies \cite{Chan_2014, Bikram_2023, Sotani_2021, sotani-kumar-UR_2021}. In this context, we focus on the universal relation between the non-radial $f$-mode frequency and the dimensionless tidal deformability for dark matter-admixed neutron stars. The top panels of Fig. Fig. \ref{fig:F-L_UR} show the normalized $f$-mode frequency ($\Bar{\omega} = \omega M$) as a function of the dimensionless tidal deformability ($\Lambda$), across three different dark matter Fermi momenta. The black dashed lines represent the best fits, and the canonical constraint from GW170817 is indicated for reference \cite{Abbott_2018}. As evident, the $f$-mode frequency remains largely independent of the selected EOS, validating the derived UR for $\Bar{\omega}$ as a function of $\Lambda$. The empirical formula used to fit the data is given by a fourth-order polynomial in $\log(\Lambda)$ as follows \cite{sotani-kumar-UR_2021,sailesh_2024,Bikram_2023}  :

\[
\Bar{\omega} = \sum_{n=0}^{4} b_n \left(\log(\Lambda)\right)^n \, ,
\]

where $b_n$ represents the fitting coefficients, which are provided in Table \ref{tab:C-L-f_coef}. The bottom panels of Fig. \ref{fig:F-L_UR} illustrate the relative deviation ($\Delta \Bar{\omega}$) between the calculated and fitted $f$-mode frequencies across the same range of $\Lambda$. The relative deviation is calculated using the formula:

\[
\Delta \Bar{\omega} = \frac{\mid{\Bar{\omega}-\Bar{\omega}_{\rm fit}}\mid}{\Bar{\omega}_{\rm fit}}.
\]

In the absence of dark matter (left panel), $\Delta \Bar{\omega}$ remains below 3\% within the canonical range of $\Lambda_{1.4}$, consistent with the GW170817 constraint. This accuracy is maintained when dark matter is introduced (middle and right panels), indicating that the presence of dark matter does not disrupt the universality of the $f$-mode frequency relation. The results suggest that the $f$-mode frequency can be reliably determined from the UR across various dark matter scenarios, with the corresponding values for different dark matter Fermi momenta listed in Table \ref{tab:C-L-f_coef}.

\subsubsection{$C-\Bar{\omega}$ relation}
The relationship between $C$ and $f$-mode frequencies was initially demonstrated by Andersson and Kokkotas \cite{Kokkotas_1999}. In this study, both with and without dark matter (DM), we use the approximation formula derived from least-squares fitting to compute the $C-\Bar{\omega}$ relations up to the fourth order as follows:

\begin{equation}\label{eq:C-F_fit}
C = \sum_{n=0}^{4} c_n \left( \Bar{\omega} \right)^n.
\end{equation}

Where $c_n$ is the fitting coefficient and the values are depicted in the Tab. \ref{tab:C-L-f_coef}.

\begin{figure*}
   \centering
   \includegraphics[width = 1.0\textwidth]{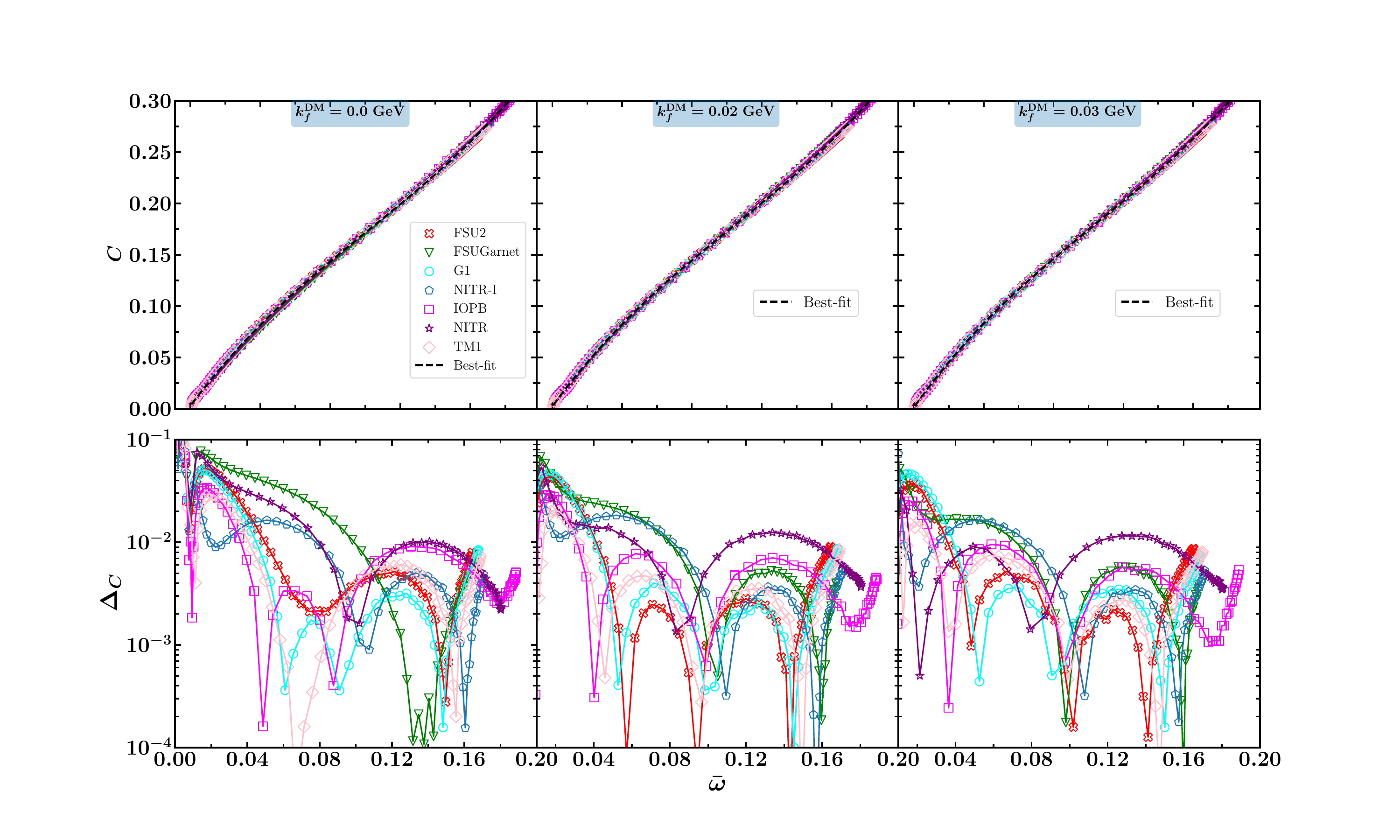}
    \caption{The $C-\Bar{\omega}$ UR is shown for both with and without DM in the upper panel. And in the lower panel, the relative deviation of compactness ($\Delta C$) is shown with normalized frequency.}
    \label{fig:C-F_UR}
\end{figure*}

In Fig. \ref{fig:C-F_UR}, we illustrate the relationship between the compactness ($C$) and the normalized $f$-mode frequency ($\Bar{\omega}$) for neutron stars, considering both scenarios with and without dark matter. The black dashed line in the top panels represents the best-fit line derived using the corresponding fitting equation \ref{eq:C-F_fit}. As shown, the UR between $C$ and $\Bar{\omega}$ holds robustly in the absence of DM, and interestingly, the inclusion of DM does not disrupt this universality, confirming the resilience of the $C-\Bar{\omega}$ relation across different dark matter scenarios. The bottom panels depict the relative deviation in compactness ($\Delta C = \frac{\mid C-C_{\text{fit}} \mid}{C_{\text{fit}}}$) with respect to $\Bar{\omega}$ for the various EOS models considered. Even with the introduction of DM, the relative deviations remain small, indicating that the compactness can still be accurately predicted from the $C-\Bar{\omega}$ UR across the tested DM Fermi momenta.

\begin{table*}
\caption{The fitting coefficients corresponding to $C-\Lambda$, $\Bar{\omega}-\Lambda$ and $C-\Bar{\omega}$ relations provided along with their reduced chi-square value for with and without DM. Also the canonical compactness from the $C-\Lambda$ UR and canonical $f$-mode frequency from $\Bar{\omega}-\Lambda$ obtained using GW170817 constraint.}
\centering
\renewcommand{\tabcolsep}{0.05cm}
\renewcommand{\arraystretch}{1.5}
\scalebox{0.65}{
    \begin{tabular}{cccccccccccc}
        \hline \hline
        \multicolumn{4}{c}{$C-\Lambda$} & \multicolumn{4}{c}{$\Bar{\omega}-\Lambda$} & \multicolumn{4}{c}{$C-\Bar{\omega}$} \\
        \hline 
        k$_f^{\rm DM}{\rm (GeV)}=$ & 0.00 & 0.03 & 0.05 &  k$_f^{\rm DM}{\rm (GeV)}=$ & 0.00 & 0.03 & 0.05 &  k$_f^{\rm DM}{\rm (GeV)}=$ & 0.00 & 0.03 & 0.05 \\
        \hline
        
        $a_0=$ & 0.350 & 0.354 & 0.358 & $b_0=$ & 0.200 & 0.200 & 0.200 & $c_0=$ & 0.002 & 0.002 & 0.001 \\
        
        $a_1\left(10^{-2}\right)=$ & -7.544 & -7.970 & -8.439 & $b_1\left(10^{-2}\right)=$ & -1.787 & -1.802 & -1.842 &$c_1=$ & 2.273 & 2.337 & 2.402 \\
        
        $a_2\left(10^{-4}\right)=$ & 1.058 & 18.248 & 31.735 & $b_2\left(10^{-2}\right)=$ & -1.636  & -1.656 & -1.708 &$c_2 =$ & -8.572 & -9.507 & -10.462 \\
        
        $a_3\left(10^{-4}\right)=$ & 9.043 & 6.364 & 4.691 & $b_3\left(10^{-3}\right)=$ & 3.746 & 3.816 & 3.991 &$c_3 =$ & 26.835 & 32.980 & 39.547 \\
        
        $a_4\left(10^{-5}\right)=$ & -5.096 & -3.725 & -2.976 & $b_4\left(10^{-4}\right)=$ & -2.306 & -2.362 & -2.497 &$c_4=$ & 5.160 & -9.616 & -26.015 \\
        
        $\chi_r^2\left(10^{-6}\right)=$ & 8.812 & 9.044 & 8.746 & $\chi_r^2\left(10^{-6}\right)=$ & 2.388 & 2.480 & 2.376 & $\chi_r^2\left(10^{-6}\right)=$ & 1.736  & 1.405 & 1.149 \\
        \hline
        \multicolumn{12}{c}{GW170817} \\
        \hline
        &  $C_{1.4}$  & $0.189^{+0.160}_{0.161}$ & $0.189^{+0.160}_{-0.161}$ & $0.188^{+0.158}_{-0.159}$ & & & $f_{1.4}({\rm kHz})$ & $2.60^{+0.456}_{-0.488}$ & $2.587^{+0.459}_{-0.490}$ & $2.549^{+0.467}_{-0.497}$  &  \\
        \hline \hline
    \end{tabular}}
    \label{tab:C-L-f_coef}
\end{table*}

\begin{figure*}
   \centering
   \includegraphics[width = 1.0\textwidth]{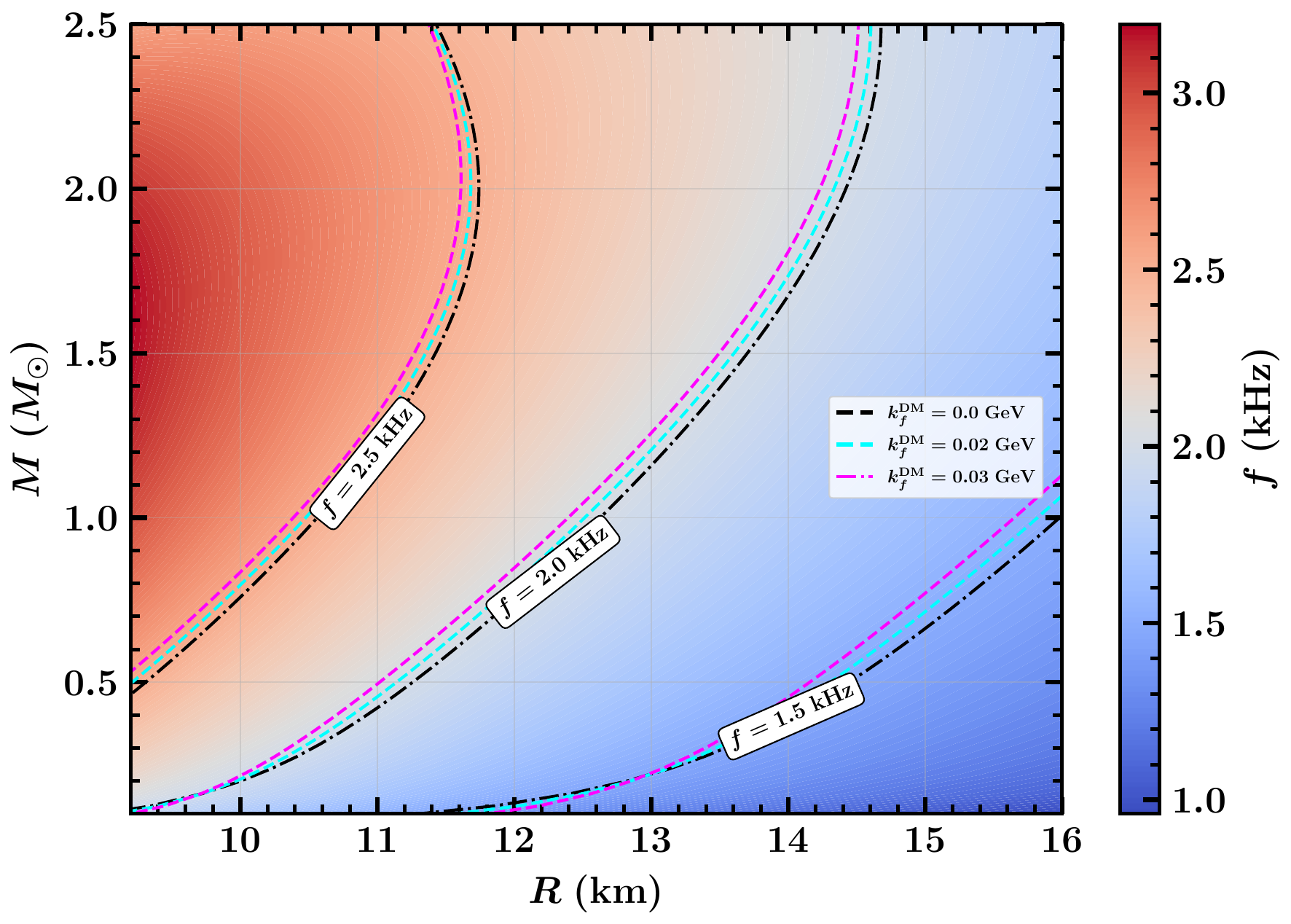}
    \caption{By employing the $C-\Bar{\omega}$ UR, the contour of frequency distribution across the $M-R$ parameter space is depicted. The distribution of frequencies is displayed by the color bar, while several dashed lines indicate variations in $k_f^{\rm DM}$.}
    \label{fig:MRContour}
\end{figure*}

By utilizing the $C-\Bar{\omega}$ UR, the distribution contour of the $f$-mode frequency over the mass-radius ($M-R$) parameter space is shown in Fig. \ref{fig:MRContour}. The color regions in the figure represent different ranges of $f$-mode frequencies, with warmer colors indicating higher frequencies and cooler colors indicating lower frequencies. The different dashed lines represent distinct combinations of mass and radius for each value of $ k_f^{\rm DM} $, predicted to exhibit the specified frequency based on the $C-\Bar{\omega}$ UR. These results highlight that in regions of low compactness, associated with cooler colors, the neutron star exhibits a lower frequency amplitude, while in regions of high compactness, associated with warmer colors, the frequency amplitude is higher. Additionally, the radius consistently decreases with the addition of DM at a constant frequency and neutron star mass, as DM softens the EOS.

\section{Summary And Conclusions}

In this study, we developed a new EDF named ``NITR-I``, which is capable of producing a maximum mass of $2.346 \, M_\odot$, consistent with the recently observed heaviest black widow pulsar, PSR J0952-0607. Our model successfully aligns with various observational constraints, including radius measurements from NICER+XMM data and the dimensionless tidal deformability constraint from the GW170817 event. We also attempted to satisfy the mass-radius constraint for the recently observed ultra-light central compact object within the supernova remnant HESS J1731-347. However, our model could not achieve this low-mass region. We extended our analysis by considering several well-established RMF models available in the literature, but none could satisfy the HESS J1731-347 constraint for a purely baryonic case.

Inspired by the study by Sagun et. al.\cite{sagun_2023hess}, we considered the possibility of DMANS to meet the HESS J1731-347 bound within a single-fluid approximation, where nucleons interact with DM candidates through the exchange of Higgs bosons, also known as non-gravitational interaction. We varied the $k_f^{\rm DM}$ parameter to explore DMANS and found that for $k_f^{\rm DM} = 0.03$ GeV, a few EOS, including our newly developed NITR-I model, successfully met the HESS J1731-347 constraints. At this DM level, the NITR-I EoS simultaneously satisfies the HESS J1731-347 constraint and the radius constraints from NICER+XMM data. However, further increases in $k_f^{\rm DM}$ cause all EoSs to satisfy the HESS bounds but fail to meet the NICER+XMM constraint. Thus, with a certain amount of DM within the NS in a single-fluid approach, we can achieve the HESS J1731-347 bounds, suggesting that this object may be a DMANS.

Additionally, based on this study and using these observational constraints, we fixed the amount of DM up to $k_f^{\rm DM} = 0.03$ GeV for further calculations. Since the presence of DM affects both the mass and radius of the NS, this conclusion led us to investigate other properties of the NS, such as tidal deformability and non-radial $f$-mode oscillations, as both are associated with the mass and radius of the NS. As the Fermi momentum of the dark matter increases, the tidal deformability decreases, while the $f$-mode oscillation frequency increases.

Furthermore, we explored various URs such as $C-\Lambda$, $\Bar{\omega}-\Lambda$, and $C-\Bar{\omega}$. We studied all the URs for both cases, with and without DM, and calculated their corresponding fitting coefficients along with the reduced chi-square values. The relative deviation for both scenarios was also calculated to estimate the accuracy of different properties. Without DM, we found that the compactness of the star was calculated with less than $3\%$ accuracy for the canonical bound of GW170817. With different DM configurations included, the universality does not break, and the estimated accuracy still remained below $3\%$ for the canonical bound. Therefore, by utilizing the GW170817 constraint for dimensionless tidal deformability, we were able to estimate the possible compactness range for the canonical NS model in our study. Similarly, we investigated non-radial $f$-mode oscillations through the $\Bar{\omega}-\Lambda$ relation and calculated the corresponding relative deviation for the $f$-mode frequency. This analysis indicated that for the GW170817 NS bound, the estimated accuracy of the $f$-mode frequency lies below $3\%$ for both cases, with and without DM. Thus, in the current study, we derived the probable range of $f$-mode frequency using the GW170817 bound for the canonical NS model. Finally, we derived the $C-\Bar{\omega}$ UR and used it to plot the distribution of $f$-mode frequency over the $M-R$ parameter space, which allowed us to identify the specific set of mass and radius for a given $f$-mode frequency.

\section{Acknowledgments}
JAP and B.K. acknowledge partial support from the Department of Science and Technology, Government of India, with grant no. CRG/2019/002691 and CRG/2021/000101 respectively. We acknowledge Vishal Parmar for his help to make ``NITR-I" EOS as unified and valuable discussions.

\bibliographystyle{ws-ijmpe}
\bibliography{hess}

\begin{thebibliography}{10}

\bibitem{pinku_jcap_2023}
P.~Routaray, S.~R. Mohanty, H.~Das, S.~Ghosh, P.~Kalita, V.~Parmar and
  B.~Kumar, {\em Journal of Cosmology and Astroparticle Physics} {\bf 2023}
  (Oct 2023)   073.

\bibitem{Abbott_2017}
 The LIGO Scientific Collaboration and Virgo Collaboration Collaborations
  (B.~P. Abbott, R.~Abbott, T.~D. Abbott {\em et~al.}), {\em Phys. Rev. Lett.}
  {\bf 119} (Oct 2017)   161101.

\bibitem{Abbott_2018}
 The LIGO Scientific Collaboration and the Virgo Collaboration Collaborations
  (B.~P. Abbott, R.~Abbott, T.~D. Abbott {\em et~al.}), {\em Phys. Rev. Lett.}
  {\bf 121} (Oct 2018)   161101.

\bibitem{Miller_2019}
M.~C. Miller, F.~K. Lamb, A.~J. Dittmann {\em et~al.}, {\em APJ} {\bf 887} (Dec
  2019)   L24.

\bibitem{Miller_2021}
M.~C. Miller, F.~K. Lamb, A.~J. Dittmann {\em et~al.}, {\em APJ} {\bf 918} (Sep
  2021)   L28.

\bibitem{RAbbott_2020}
R.~Abbott, T.~D. Abbott, S.~Abraham {\em et~al.}, {\em The Astrophysical
  Journal} {\bf 896} (Jun 2020)   L44.

\bibitem{Riley-nicer_2019}
T.~E. Riley, A.~L. Watts, S.~Bogdanov, P.~S. Ray, R.~M. Ludlam, S.~Guillot,
  Z.~Arzoumanian, C.~L. Baker, A.~V. Bilous, D.~Chakrabarty, K.~C. Gendreau,
  A.~K. Harding, W.~C.~G. Ho, J.~M. Lattimer, S.~M. Morsink and T.~E.
  Strohmayer, {\em The Astrophysical Journal Letters} {\bf 887} (Dec 2019)
  L21.

\bibitem{Miller-nicer_2019}
M.~C. Miller, F.~K. Lamb, A.~J. Dittmann, S.~Bogdanov, Z.~Arzoumanian, K.~C.
  Gendreau, S.~Guillot, A.~K. Harding, W.~C.~G. Ho, J.~M. Lattimer, R.~M.
  Ludlam, S.~Mahmoodifar, S.~M. Morsink, P.~S. Ray, T.~E. Strohmayer, K.~S.
  Wood, T.~Enoto, R.~Foster, T.~Okajima, G.~Prigozhin and Y.~Soong, {\em The
  Astrophysical Journal Letters} {\bf 887} (Dec 2019)   L24.

\bibitem{Romani_2022}
R.~W. Romani, D.~Kandel, A.~V. Filippenko, T.~G. Brink and W.~Zheng, {\em The
  Astrophysical Journal Letters} {\bf 934} (Jul 2022)   L17.

\bibitem{HESS_2022}
V.~{Doroshenko}, V.~{Suleimanov}, G.~{P{\"u}hlhofer} and A.~{Santangelo}, {\em
  Nature Astronomy} {\bf 6} (December 2022) 1444.

\bibitem{Horvath-hess_2023}
J.~E. Horvath, L.~S. Rocha, L.~M. de~S{\'{a} }, P.~H. R.~S. Moraes, L.~G.
  Bar{\~{a}}o, M.~G.~B. de~Avellar, A.~Bernardo and R.~R.~A. Bachega, {\em
  Astronomy \& Astrophysics} {\bf 672} (Apr 2023)   L11.

\bibitem{clemente_2023hess}
F.~D. Clemente, A.~Drago and G.~Pagliara, {\em The Astrophysical Journal} {\bf
  967} (May 2024)   159.

\bibitem{hcdas_2023hess}
H.~C. Das and L.~L. Lopes, {\em Mon. Not. R. Astron. Soc.} {\bf 525} (November
  2023) 3571.

\bibitem{rather_2023quark}
I.~A. Rather, G.~Panotopoulos and I.~Lopes, {\em Eur. Phys. J. C} {\bf 83}
  (November 2023) 1.

\bibitem{kurbis_2019_prcmodel}
N.~Zabari, S.~Kubis and W.~W\'ojcik, {\em Phys. Rev. C} {\bf 100} (Jul 2019)
  015808.

\bibitem{kubis_2023hess}
S.~Kubis, W.~W\'ojcik, D.~A. Castillo and N.~Zabari, {\em Phys. Rev. C} {\bf
  108} (Oct 2023)   045803.

\bibitem{huang_2023hess}
K.~Huang, H.~Shen, J.~Hu and Y.~Zhang, {\em Phys. Rev. D} {\bf 109} (Feb 2024)
   043036.

\bibitem{sagun_2023hess}
V.~Sagun, E.~Giangrandi, T.~Dietrich, O.~Ivanytskyi, R.~Negreiros and
  C.~Providência, {\em The Astrophysical Journal} {\bf 958} (Nov 2023)  ~49.

\bibitem{McCullough_WD-DM_2010}
M.~McCullough and M.~Fairbairn, {\em Phys. Rev. D} {\bf 81} (Apr 2010)
  083520.

\bibitem{Goldman_1989}
I.~Goldman and S.~Nussinov, {\em Phys. Rev. D} {\bf 40} (Nov 1989) 3221.

\bibitem{Sandin_2009}
F.~Sandin and P.~Ciarcelluti, {\em Astroparticle Phys.} {\bf 32} (Dec 2009)
  278–284.

\bibitem{Kouvaris_2008}
C.~Kouvaris, {\em Phys. Rev. D} {\bf 77} (Jan 2008)   023006.

\bibitem{Kouvaris_2011}
C.~Kouvaris and P.~Tinyakov, {\em Phys. Rev. D} {\bf 83} (Apr 2011)   083512.

\bibitem{Grigorious_2017}
G.~Panotopoulos and I.~Lopes, {\em Phys. Rev. D} {\bf 96} (Oct 2017)   083004.

\bibitem{Lopes-harish_2023}
L.~L. Lopes and H.~Das, {\em Journal of Cosmology and Astroparticle Physics}
  {\bf 2023} (May 2023)   034.

\bibitem{Kokkotas_2001}
K.~D. Kokkotas and J.~Ruoff, {\em Astronomy \& Astrophysics} {\bf 366} (Feb
  2001) 565.

\bibitem{souhardya_2023}
S.~Sen, S.~Kumar, A.~Kunjipurayil, P.~Routaray, S.~Ghosh, P.~J. Kalita, T.~Zhao
  and B.~Kumar, {\em Galaxies} {\bf 11}  (2023).

\bibitem{pinku_prd_2023}
P.~Routaray, H.~C. Das, S.~Sen, B.~Kumar, G.~Panotopoulos and T.~Zhao, {\em
  Phys. Rev. D} {\bf 107} (May 2023)   103039.

\bibitem{Ishfaq_rad-delta_2023}
I.~A. Rather, K.~D. Marquez, G.~Panotopoulos and I.~Lopes, {\em Phys. Rev. D}
  {\bf 107} (Jun 2023)   123022.

\bibitem{Bikram_2021-fmode}
B.~K. Pradhan and D.~Chatterjee, {\em Phys. Rev. C} {\bf 103} (Mar 2021)
  035810.

\bibitem{athul_2022}
A.~Kunjipurayil, T.~Zhao, B.~Kumar, B.~K. Agrawal and M.~Prakash, {\em Phys.
  Rev. D} {\bf 106} (Sep 2022)   063005.

\bibitem{sailesh_2024}
S.~R. Mohanty, S.~Ghosh, P.~Routaray, H.~Das and B.~Kumar, {\em Journal of
  Cosmology and Astroparticle Physics} {\bf 2024} (Mar 2024)   054.

\bibitem{Probit_2024}
P.~J. Kalita, P.~Routaray, S.~Ghosh, B.~Kumar and B.~Agrawal, {\em Journal of
  Cosmology and Astroparticle Physics} {\bf 2024} (Apr 2024)   065.

\bibitem{Kent_yagi_2013}
K.~Yagi and N.~Yunes, {\em Phys. Rev. D} {\bf 88} (Jul 2013)   023009.

\bibitem{Kent_yagi_2015}
K.~Yagi and N.~Yunes, {\em Phys. Rev. D} {\bf 91} (Jun 2015)   123008.

\bibitem{Gupta_2017}
T.~Gupta, B.~Majumder, K.~Yagi {\em et~al.}, {\em Classical and Quantum
  Gravity} {\bf 35} (Dec 2017)   025009.

\bibitem{Yeung_2021}
C.-H. Yeung, L.-M. Lin, N.~Andersson {\em et~al.}, {\em Universe} {\bf 7}
  (2021).

\bibitem{Harish__I_LOVE_C_2022}
H.~C. Das, {\em Phys. Rev. D} {\bf 106} (Nov 2022)   103518.

\bibitem{sotani_2011-cowling}
H.~Sotani, N.~Yasutake, T.~Maruyama and T.~Tatsumi, {\em Phys. Rev. D} {\bf 83}
  (Jan 2011)   024014.

\bibitem{Flores_2014-cowling}
C.~V. Flores and G.~Lugones, {\em Classical and Quantum Gravity} {\bf 31} (Jul
  2014)   155002.

\bibitem{FURNSTAHL_1996}
R.~Furnstahl, B.~D. Serot and H.-B. Tang, {\em Nuclear Physics A} {\bf 598}
  (1996) 539.

\bibitem{Frun_1997}
R.~J. Furnstahl, B.~D. Serot and H.-B. Tang, {\em Nucl. Phys. A} {\bf 615}
  (1997) 441.

\bibitem{singh_2014}
S.~K. Singh, S.~K. Biswal, M.~Bhuyan and S.~K. Patra, {\em Journal of Physics
  G: Nuclear and Particle Physics} {\bf 41} (Mar 2014)   055201.

\bibitem{Kumar_2017}
B.~Kumar, S.~Singh, B.~Agrawal and S.~Patra, {\em Nuclear Physics A} {\bf 966}
  (2017) 197.

\bibitem{Kumar_2018}
B.~Kumar, S.~K. Patra and B.~K. Agrawal, {\em Phys. Rev. C} {\bf 97} (Apr 2018)
    045806.

\bibitem{Parmar_2022_1}
V.~Parmar, H.~C. Das, A.~Kumar, M.~K. Sharma and S.~K. Patra, {\em Phys. Rev.
  D} {\bf 105} (Feb 2022)   043017.

\bibitem{BKAgrawal_2005}
B.~K. Agrawal, S.~Shlomo and V.~K. Au, {\em Phys. Rev. C} {\bf 72} (Jul 2005)
  014310.

\bibitem{BKAgrawal_2006}
R.~Kumar, B.~K. Agrawal and S.~K. Dhiman, {\em Phys. Rev. C} {\bf 74} (Sep
  2006)   034323.

\bibitem{Parmar_2022}
V.~Parmar, H.~C. Das, A.~Kumar, M.~K. Sharma, P.~Arumugam and S.~K. Patra, {\em
  Phys. Rev. D} {\bf 106}  (2022)   023031.

\bibitem{Drischler_2021}
C.~Drischler, S.~Han, J.~M. Lattimer, M.~Prakash, S.~Reddy and T.~Zhao, {\em
  Phys. Rev. C} {\bf 103} (Apr 2021)   045808.

\bibitem{arpan_2019}
A.~Das, T.~Malik and A.~C. Nayak, {\em Phys. Rev. D} {\bf 99} (Feb 2019)
  043016.

\bibitem{harishmnras_2020}
H.~C. Das, A.~Kumar, B.~Kumar, S.~K. Biswal, T.~Nakatsukasa, A.~Li and S.~K.
  Patra, {\em Monthly Notices of the Royal Astronomical Society} {\bf 495} (05
  2020) 4893.

\bibitem{pinku_mnras_2023}
P.~Routaray, A.~Quddus, K.~Chakravarti and B.~Kumar, {\em Monthly Notices of
  the Royal Astronomical Society} {\bf 525} (09 2023) 5492.

\bibitem{Das_fmode_2021}
H.~C. Das, A.~Kumar, S.~K. Biswal and S.~K. Patra, {\em Phys. Rev. D} {\bf 104}
  (Dec 2021)   123006.

\bibitem{Tolman_1939}
R.~C. Tolman, {\em Phys. Rev.} {\bf 55} (Feb 1939) 364.

\bibitem{Oppenheimer_1939}
J.~R. Oppenheimer and G.~M. Volkoff, {\em Phys. Rev.} {\bf 55} (Feb 1939) 374.

\bibitem{Hinderer_2008}
T.~Hinderer, {\em The Astrophysical Journal} {\bf 677} (Apr 2008) 1216.

\bibitem{Kumartidal_2017}
B.~Kumar, S.~K. Biswal and S.~K. Patra, {\em Phys. Rev. C} {\bf 95} (Jan 2017)
   015801.

\bibitem{Chen_2014}
W.-C. Chen and J.~Piekarewicz, {\em Phys. Rev. C} {\bf 90} (Oct 2014)   044305.

\bibitem{Chen_2015}
W.-C. Chen and J.~Piekarewicz, {\em Physics Letters B} {\bf 748}  (2015) 284.

\bibitem{Sugahara_1994}
Y.~Sugahara and H.~Toki, {\em Nucl. Phys. A} {\bf 579} (October 1994) 557.

\bibitem{Cromartie_2020}
H.~T. Cromartie, E.~Fonseca, S.~M. Ransom {\em et~al.}, {\em Nature Astronomy}
  {\bf 4} (Jan 2020) 72.

\bibitem{Andersson-Kokkotas_GW_1996}
N.~Andersson and K.~D. Kokkotas, {\em Phys. Rev. Lett.} {\bf 77} (Nov 1996)
  4134.

\bibitem{sotani-kumar-UR_2021}
H.~Sotani and B.~Kumar, {\em Phys. Rev. D} {\bf 104} (Dec 2021)   123002.

\bibitem{Andersson-Kokkotas_asteroseismology_1998}
N.~Andersson and K.~D. Kokkotas, {\em Monthly Notices of the Royal Astronomical
  Society} {\bf 299} (10 1998) 1059.

\bibitem{Chakrabarti_2014}
S.~Chakrabarti, T.~Delsate, N.~G\"urlebeck {\em et~al.}, {\em Phys. Rev. Lett.}
  {\bf 112} (May 2014)   201102.

\bibitem{Haskell_2013}
B.~Haskell, R.~Ciolfi, F.~Pannarale {\em et~al.}, {\em Monthly Notices of the
  Royal Astronomical Society: Letters} {\bf 438} (12 2013) L71.

\bibitem{Bandyopadhyay2018}
D.~Bandyopadhyay, S.~A. Bhat, P.~Char {\em et~al.}, {\em The European Physical
  Journal A} {\bf 54} (Feb 2018)  ~26.

\bibitem{Prashant_2024}
P.~Thakur, T.~Malik, A.~Das, T.~K. Jha and C.~m.~c. Provid\^encia, {\em Phys.
  Rev. D} {\bf 109} (Feb 2024)   043030.

\bibitem{Maselli_C-L_2013}
A.~Maselli, V.~Cardoso, V.~Ferrari, L.~Gualtieri and P.~Pani, {\em Phys. Rev.
  D} {\bf 88} (Jul 2013)   023007.

\bibitem{Chan_2014}
T.~K. Chan, Y.-H. Sham, P.~T. Leung {\em et~al.}, {\em Phys. Rev. D} {\bf 90}
  (Dec 2014)   124023.

\bibitem{Bikram_2023}
B.~K. Pradhan, A.~Vijaykumar and D.~Chatterjee, {\em Phys. Rev. D} {\bf 107}
  (Jan 2023)   023010.

\bibitem{Sotani_2021}
H.~Sotani and B.~Kumar, {\em Phys. Rev. D} {\bf 104} (Dec 2021)   123002.

\bibitem{Kokkotas_1999}
K.~D. Kokkotas and B.~G. Schmidt, {\em Living Reviews in Relativity} {\bf 2}
  (Sep 1999)  ~2.

\end{thebibliography}

\end{document}